\documentclass[12pt, twoside, leqno]{article}
\usepackage{latexsym}
\usepackage{amsmath}
\usepackage{amsfonts}
\usepackage{amssymb}
\usepackage{amsthm}
\usepackage{color}
\numberwithin{equation}{section}

\input xy
\xyoption{all} \hfuzz=1.5pt \tolerance=500
\theoremstyle{plain}

\theoremstyle{definition}

\newcommand{\tr}{\triangleright}
\newcommand{\tl}{\triangleleft}
\newcommand{\cd}{{\cdot}}
\newcommand{\ot}{{\otimes}}

\newcommand{\cK}{{\cal F}}
\newcommand{\la}{\langle}
\newcommand{\ra}{\rangle}
\newcommand{\id}{{\bf 1}}
\newcommand{\di}{\diamondsuit}
\newcommand{\q}{\quad}
\newcommand{\qq}{\qquad}
\newcommand{\va}{\varphi}
\newcommand{\rr}{{\cal R}}
\newcommand{\ii}{\infty}
\newcommand{\mt}{\mapsto}
\newcommand{\co}{{\mathbb C}}

\newcommand{\bb}{{\cal B}}

\newcommand{\al}{{\alpha}}     
     
\newcommand{\de}{{\delta}}     \newcommand{\lm}{{\lambda}}

\newcommand{\N}{{\mathbb N}}  
  
\def\x{\relax\ifmmode {\mbox{*}}\else*\fi}

\newcommand{\ed}{\end{document}}

\begin{document}

\bigskip
\centerline{{\bf A.~Ya.~Helemskii.  }}

%\pagestyle{myheadings} \markboth{\centerline{\small{\sc
%            A.~Ya.~Helemskii}}}
%         {\centerline{\small{\sc Structures on the way  from classical to quantum spaces}}}
 %\title{\bf 
\bigskip
\centerline{{\bf{\large Multi-normed spaces, based on non-discrete}}}

\centerline{{\bf{\large  measures, and their tensor products}}\footnote{The research was supported 
by the Russian Foundation for Basic Researches (grunt No. 15-01-08392)}} 

\bigskip
\begin{abstract}

\medskip 
%{\small 
It was A.~Lambert who discovered a new type of structures, situated, in a sense, between normed 
spaces and (abstract) operator spaces. His definition was based on the notion of amplification a 
normed space by means of spaces $\ell_2^n$. Afterwards several mathematicians investigated more 
general structure, ``$p$--multi-normed space'', introduced with the help of spaces $\ell_p^n;1\le 
p\le\ii$. In the present paper we pass from $\ell_p$ to $L_p(X,\mu)$ with an arbitrary measure. 
This happened to be possible in the frame-work of the non-coordinate (``index-free'') approach to 
the notion of amplification, equivalent in the case of a discrete counting measure to the approach 
in mentioned articles. 

Two categories arise. One consists of amplifications by means of an arbitrary normed space, and 
another one consists $p$--convex amplifications by means of $L_p(X,\mu)$. Each of them has its own 
tensor product of its objects whose existence is proved by a respective explicit construction. As a 
final result, we show that the ``$p$--convex'' tensor product has especially transparent form for 
the so-called minimal $L_p$--amplifications of $L_q$--spaces, where $q$ is the conjugate of $p$. 
Namely, tensoring $L_q(Y,\nu)$ and $L_q(Z,\lm)$, we get $L_q(Y\times Z,\nu\times\lm)$. 

%Библиография: 29 названий.

\end{abstract}



\section{Introduction}

The subject of the present paper is a rather far-reaching (through several steps) generalization  
of the structure, introduced in the PhD thesis of A.~Lambert~\cite{lam}; his superviser was 
 G.~Wittstock, one of the founding fathers of operator space theory. (See also~\cite{loi}). 
 This structure is, speaking informally, intermediate between the classical structure of a normed space
 and the structure of an abstract operator space. 
 %or, as they also say, of a quantum space. 
 (The latter
 is presented in widely known textbooks~\cite{efr,pis,paul,blem}; see also~\cite{heb2}). 
 Lambert suggested to endow a given linear space $E$ by a sequence of norms on 
  spaces $E^n;n\in\N$, that is on  columns of all 
 possible sizes, consisting of vectors from $E$. 
 % (That is, not of norms on matrices, as in the theory of operator spaces but norms of 
 %their first columns). 
 (Thus, he deals with norms on columns, and not on matrices, as in the theory of operator spaces).
  These norms must satisfy two axioms, 
 ``contractiveness'' and ``convexity'' that were formulated in terms of the spaces $\ell_2^n;n\in{\N}$. 
Lambert called the resulting objects ``Operatorfolgenr\"aume''. In the respective rising category 
Lambert has constructed two tensor products, ``maximal''~\cite[3.1.1]{lam} and 
``minimal''~\cite[3.1.3]{lam}; the former one can be considered as a predecessor of the tensor 
product, introduced in Section 5 of the present paper. 

The theory of Lambert had various connections with the classical theory of normed spaces as well as 
with the theory of operator spaces, shedding in many occasions a new light in their relationship.  
Later a team of mathematicians, embarking from essentially different problems, related to Banach 
lattices, came to more general structures. However, it was done again in the frame-work of the 
``coordinate approach'', based on the consideration of columns of arbitrary size. 
%.
%-wise amplifications'' of a given linear space, that is on columns of arbitrary 
% sizes, consisting of its vectors. 
First, there were Dales and Polyakov~\cite{dpol}, 
 soon after joined by Daws, Pham, Ramsden, Laustsen, Oikhberg and Troitsky~\cite{ddaws, ddaws2,dal}. 
These authors created rich and ramified theory, from which we are most interesting in the so-called 
  $p$--multi-normed spaces~\cite{dal}: those satisfying the analogue of the contractiveness axiom of 
  Lambert,
but now in terms of the spaces $\ell_p^n$ with arbitrary fixed $p;1\le p\le\ii$. The ``best'' of 
these structures satisfy also the analogue of the convexity axiom, the so-called $p$--convexity. 
(Lambert has $p=2$). 

The present paper pursuers two aims. First, we extend the class of the structures in question, 
passing from $p$-multi-normed spaces to their ``continuous'' (non-discrete) versions. Namely, we 
%pass, in the capacity of a base space, from $\ell_p$ to the spaces $L_p(X,\mu)$ with arbitrary 
change, in the capacity of a base space, $\ell_p$ to $L_p(X,\mu)$ with arbitrary measure. This 
becomes possible, if we replace the coordinate approach by the so-called non-coordinate approach to 
what we call an amplification. 
%pass from the coordinate approach to the so-called non-coordinate approach to what we call an amplification. 

In the context of operator spaces the latter approach was known to specialists, and it was 
systematically presented in~\cite{heb2} (see also~\cite{he2}). In the context of Lambert spaces it 
was applied in~\cite{heh}, where several notions and facts from the present paper have their 
prototypes. The essence of this approach is as follows. Instead of a sequence of norms on all $E^n$ 
%columns of vectors from $E$, 
%, each one on its own column space, on column spaces, consisting of vectors from $Е$, 
we consider a norm on a single space ${\bf L}\ot E$, where ${\bf L}$ is a chosen and fixed ``base'' 
space.
%, e.g., the same $L_p(X,\mu)$. 
Such an approach in the case, where ${\bf L}:=L_p(X,\mu)$ with discrete counting measure, is 
equivalent to the coordinate approach, accepted in the above-cited papers. However, it seems that 
as a whole it provides greater possibilities. In the frame-work of this approach the axiom of 
contractiveness transforms to the condition on the normed space ${\bf L}\ot E$ to be a contractive 
module over $\bb({\bf L})$. As to the ``non-coordinate'' version of the axiom of $p$-convexity, it 
can be defined under certain assumptions on ${\bf L}\ot E$ that make this space very similar to 
$L_p(X,\mu)$. 

Spaces of the form ${\bf L}\ot E$ are called in this paper  ${\bf L}$--spaces. Most of all, we are 
interested in their tensor products.
% of arising objects, the latter called in this paper $\bf L}$--spaces. 
We introduce two essentially different varieties of this notion. The first one, 
 `` $\ot_{\bf L}$ '', is defined for the case of general base spaces ${\bf L}$ that are endowed with a 
certain additional structure, a bilinear operator $\di:{\bf L}\times{\bf L}\to{\bf L}$, possessing 
some natural properties. Another kind of tensor product, denoted by `` $\ot_{p{\bf L}}$ '' and 
called $p$--convex, is constructed 
% only
 for the class of $p$--convex ${\bf L}$-spaces and only in the case, when 
our base space is $L_p(X,\mu)$. Each of these tensor products is defined in terms of universal 
property for the respective class of bilinear operators, and its existence theorem is proved by 
displaying its own explicit construction. We present several examples. In particular, we show that 
 the tensor product ``$\ot_{\bf L}$ '' acquires sufficiently transparent guise for ${\bf L}$--spaces
that have the biggest of all possible norms. What is, perhaps, more interesting is that the second 
tensor product acquires a transparent concrete guise for $L_p(X,\mu)$--spaces of the form 
$L_q(\cd)$, this time endowed 
 with the minimal norm; here $1<p<\ii$ and $q$ is the conjugate to $p$. Namely, up to a 
 $L_p(X,\mu)$--isometric 
isomorphism of $L_p(X,\mu)$-пространств (which is an analogue of complete isometric isomorphism of 
operator spaces), we have 
\[
L_q(Y,\nu)\ot_{p{\bf L}}L_q(Z,\lm)=L(Y\times Z,\nu\times\lm). 
\]

The author is much indebted to N.T.Nemesh for valuable discussions. 

\section{${\bf L}$--spaces and ${\bf L}$--quantization. }

As usual, we denote by $\bb(E,F)$ the space of all bounded operators between the normed spaces $E$ 
and $F$, and consider it with the operator norm. We write $\bb(E)$ instead of $\bb(E,E)$. The 
identity operator on $E$ will be denoted by $\id_E$. Two projections $P$ and $Q$ on $E$ are called 
{\it orthogonal}, if $PQ=QP=0$. 

The symbol $\ot$ is used for the algebraic tensor product of linear spaces and for 
 elementary tensors. The symbols $\ot_{pr}$ and $\ot_{in}$ denote the non--completed projective and injective 
 tensor product of normed spaces, respectively.
 
Choose and fix (so far arbitrary) normed space ${\bf L}$, which we shall call the {\it bаse space}. 
Let us write $\bb$ instead of $\bb({\bf L})$. 

\medskip
In what follows we need the triple notion of the so-called amplification. First, we amplify linear 
spaces, then linear operators and finally bilinear operators. Note that these amplifications differ 
from the amplifications from~\cite{heb2}, serving in the theory of operator spaces. 
%(cf.~\cite{heb2}). 

The {\it amplification} of a given linear space $E$ is the tensor product ${\bf L}\ot E$. Usually 
we briefly denote it by ${\bf L}E$, and an elementary tensor, say $\xi\ot x; \xi\in{\bf L}, x\in 
E$, by $\xi x$. Note that ${\bf L}E$ is a left module over the algebra $\bb$ with the outer 
multiplication `` $\cd$ '', well defined by $a\cd(\xi x):=a(\xi)x$. 

\medskip
{\bf Definition 2.1.} A semi-norm on ${\bf L}E$ is called {\it ${\bf L}$--seminorm on} $E$, if the 
left $\bb$-module ${\bf L}E$ is contractive, that is if we always have the estimate $\|a\cd 
u\|\le\|a\|\|u\|$. The space $E$, endowed by a ${\bf L}$--seminorm, is called {\it seminormed ${\bf 
L}$--space}). If the seminorm in question is actually a norm, we speak, naturally, about a {\it 
normed ${\bf L}$--space}, and in this case we usually omit the word ``normed''. 

\medskip
{\bf Example 2.2}. Let $(X,\mu)$ be a measure space. Our principal example of a base space is 
$L_p(X,\mu)$, where $1\le p\le\ii$. 
%; мы всегда будем считать, что наша мера такова, что это пространство бесконечномерно. 
(As the main reference in the measure theory, we shall use the textbook~\cite{bog}). For simplicity 
{\bf we always assume that all our measures have a countable basis.} (Thus, the set of atoms is no 
more than countable). If there is no danger of misunderstanding, we shall speak about a measure 
space $X$ and the normed space $L_p(X)$; in particular, $X\times Y$ denotes the cartesian product 
of respective measure spaces.  

\medskip
{\bf Remark 2.3.} As to the former papers, cited above, they consider, after translation into the 
``index-free'' language, the case $X:=\N$ with the counting measure. In particular, spaces of 
Lambert are those with ${\bf L}=\ell_2$, whereas in spaces of Dales/Polyakov we have 
 ${\bf L}=\ell_\ii$ or ${\bf L}=\ell_1$. Finally, if ${\bf L}=\ell_p$, then the notion of ${\bf 
 L}$--space is equivalent to the notion of $p$--multi-normed space of 
 Dales/Laustsen/Oikberg/Troitsky~\cite[2.2]{dal} 
%в пространствах Дэйлса--Лаустсена--Ойкберга--Троицкого ${\bf L}=\ell_p$ при произвольных $p$. 

\medskip
A seminormed ${\bf L}$--space $E$ becomes seminormed space in the usual sense, if  for $x\in E$ we 
set $\|x\|:=\|\xi x\|$, where $\xi\in{\bf L}$ is an arbitrary vector with $\|\xi\|=1$. Clearly, the 
result does not depend on a choice of $\xi$. The obtained seminormed space is called {\it 
underlying space} of a given ${\bf L}$-space, and the latter is called an {\it ${\bf 
L}$--quantization} of a former. (We use such a term by analogy with quantizations in operator space 
theory; see, e.g.,~\cite{ef5},~\cite{efr} or~\cite{heb2}). Obviously, for all $\xi\in{\bf L}$ and 
$x\in E$ we have $\|\xi x\|=\|\xi\|\|x\|$. 

\medskip
It is easy to verify that the space of scalars, $\co$, has the only ${\bf L}$--quantization, given 
by the identification of ${\bf L}\co$ with ${\bf L}$. 

\medskip
{\bf Proposition 2.4}. {Let $E$ be a seminormed ${\bf L}$--space with a normed underlying space. 
Then the ${\bf L}$--seminorm on ${\bf L}$ is a norm.} 

\smallskip
$\tl$ Take $u\in{\bf L}E; u\ne0$ and represent it as $\sum_{k=1}^n\xi_kx_k$, where $\xi_k$ are 
linearly independent, $\|\xi_1\|=1$ and $x_1\ne0$. Obviously, there exists  $T\in\bb$ with 
$\|T\xi_1\|=1$ and $T\xi_k=0$ for $k>1$. Then, according to Definition 2.1, we have  
$\|T\|\|u\|\ge\|T\cd u\|=\|x_1\|>0. \q \tr$ 
 
\medskip
{\bf Example 2.5.}  Every normed space, say $E$, has, generally speaking, a lot of ${\bf 
L}$--quantizations. We distinguish two of them. The 
 ${\bf L}$--space, denoted by $E_{\max}$, respectively $E_{\min}$, has the ${\bf L}$--norm, 
 obtained by the endowing ${\bf L}E$ with the norm of ${\bf L}\ot_{pr}E$,
respectively of ${\bf L}\ot_{in}E$. We denote the norm on the former and on the latter space by 
$\|\cd\|_{\max}$ and $\|\cd\|_{\min}$, respectively; accordingly, the corresponding ${\bf 
L}$-quantizations of $E$ will be called {\it maximal} and {\it minimal}. Clearly, the ${\bf 
L}$-norm of $E_{\max}$ is the greatest of all ${\bf L}$-norms of ${\bf L}$-quantizations of $E$. 
The adjective ``minimal'' will be justified a little bit later. 

\medskip
{\bf Example 2.6.} We want to introduce an ${\bf L}$--quantization of the ``classical'' projective 
tensor product $E\ot_{pr} F$ of two normed spaces, when one of tensor factors, say, to be definite, 
$F$, is itself an ${\bf L}$--space. 

Consider the linear isomorphism $\beta: {\bf L}(E\ot F)\to E\ot_{pr}({\bf L}F):\xi(x\ot y)\mt 
x\ot\xi y$ and introduce a norm on ${\bf L}(E\ot F)$ by setting $\|U\|':=\|\beta(U)\|$. The space 
$E\ot_{pr}({\bf L}F)$, being a projective tensor product of a normed space and a contractive 
$\bb$-module, has itself a standard structure of a contractive $\bb$-module. Since $\beta$ is a 
$\bb$-module morphism, the same is true with ${\bf L}(E\ot F)$. Thus $E\ot F$ becomes an ${\bf 
L}$--space. Denote the norm of the respective underlying space just by $\|\cd\|$, and the norm on в 
$E\ot_{pr}F$ by $\|\cd\|_{pr}$. We must show that $\|\cd\|=\|\cd\|_{pr}$. 

Take an arbitrary $u\in E\ot F$. Since it is clear that $\|\cd\|$  is a cross-norm, we have  
$\|u\|\le\|u\|_{pr}$. It remains to show that for every $\xi\in{\bf L},\|\xi\|=1$ we have $\|\xi 
u\|'\ge\|u\|_{pr}$. 

Identifying $\bb$-modules ${\bf L}(E\ot F)$ and $E\ot_{pr}({\bf L}F)$ by means of $\beta$, we
 represent $\xi u$ as $\sum_{k=1}^nx_k\ot w_k; x_k\in E,w_k\in{\bf L}F$. Take a functional
$f:{\bf L}\to\co$ such that $f(\xi)=1, \|f\|=1$, and the operator $T:{\bf L}\to{\bf L}:\eta\mt 
f(\eta)\xi$; clearly, $\|T\|=1$. It is obvious that $T\cd w_k=\xi y_k$ for some $y_k\in F; 
k=1,...,n$. Therefore we have $\sum_{k=1}^n\|x_k\|\|w_k\|\ge \sum_{k=1}^n\|x_k\|\|T\cd w_k\|= 
\sum_{k=1}^n\|x_k\|\|y_k\|$. But  $\xi u=T\cd(\xi u)=\sum_{k=1}^nx_k\ot T\cd w_k= 
\xi(\sum_{k=1}^nx_k\ot y_k)$. Consequently,  $u=\sum_{k=1}^nx_k\ot y_k$. It follows that 
$\sum_{k=1}^n\|x_k\|\|w_k\|\ge\|u\|_{pr}$, and this, because of the definition of the projective 
norm, implies the desired estimate $\|\xi u\|'\ge\|u\|_{pr}$. 

\medskip
{\bf Remark 2.7.} In the textbook~\cite{heb2} the amplification of a space $E$ is defined as ${\cal 
F}\ot E$, where $\cK$ is the space of bounded finite rank operators on a certain Hilbert space $L$. 
However, it is worthy to mention that the norm on ${\cal F}\ot E$, making $E$ an (abstract) 
operator space, is not always a $\bf L$-norm in the sense of Definition 2.1. The simplest 
counter-example is $E:=\bb(L)$ with the standard norm of an operator space. Using the estimates, 
obtained by Tomiyama~\cite{tom}, one can prove the following assertion: if $\tau:\cK\to\cK$ is the 
operator, acting as the transpose of corresponding matrices, then for every $C>0$ there exists 
$u\in\cK E$ such that, notwithstanding $\|\tau\|=1$, we have $\|\tau\cd u\|>C\|u\|$. 

\section{${\bf L}$--bounded linear and bilinear operators}

Suppose we are given an operator $\va:E\to F$ between linear spaces. Denote, for brevity, the 
operator $\id\ot\va:{\bf L}E\to{\bf L}F$ (taking $\xi x$ to $\xi\va(x)$) by $\va_\ii$ and call it 
{\it amplification} of $\va$. Obviously, $\va_\ii$ is a morphism of left $\bb$-modules. 

\medskip
{\bf Definition 3.1.} An operator $\va:E\to F$ between seminormed ${\bf L}$--spaces is called {\it 
${\bf L}$--bounded}, if the operator $\va_\ii$ is bounded. Then we set $\|\va\|_{{\bf 
L}b}:=\|\va_\ii\|$. Similarly, in terms of $\va_\ii$ we define {\it ${\bf L}$--contractive} and 
{\it ${\bf L}$--isometric} operator, and also {\it ${\bf L}$--isometric isomorphism}. 

\medskip
If $\va$ is bounded, being considered between the respective underlying seminormed spaces, we say 
that it is (just) {\it bounded}, and denote its operator seminorm, as usual, by $\|\va\|$. Clearly, 
every $L$--bounded operator $\va:E\to F$ is bounded, and $\|\va\|\le\|\va\|_{{\bf L}b}$ . 

\medskip
Some operators between ${\bf L}$--spaces, bounded as operators between underlying spaces, are 
``automatically'' ${\bf L}$--bounded. Here is the first phenomenon of that kind. 

\medskip
{\bf Proposition 3.2.} {\it Let $E$ be a ${\bf L}$--space. Then every bounded functional 
$f:E\to\co$ is ${\bf L}$--bounded, and $\|f\|_{{\bf L}b}:=\|f\|$. } 

\smallskip 
$\tl$ Recall that we have $\|f_\ii(u)\|=\max\{|\la f_\ii(u),\al\ra|; \al\in{\bf L}^*, \|\al\|=1\}$. 
Set, for every $\al$, $x_\al:=(\al\ot\id_E)(u)\in E$; then, representing $u$ as a sum of elementary 
tensors, we obtain that $\la f_\ii(u),\al\ra=f(x_\al)$. Now fix an arbitrary $\eta\in{\bf 
L};\|\eta\|=1$ and consider the operator $S:{\bf L}\to{\bf L}:\xi\mt\al(\xi)\eta$; obviously, we 
have $\|S\|=\|\al\|=1$. But the same representation of $u$ implies $S\cd u=\eta x_\al$. Therefore 
we have $\|x_\al\|=\|\eta x_\al\|\le\|S\|\|u\|=\|u\|$. It follows that $|\la 
f_\ii(u),\al\ra|\le\|f\|\|x_\al\|=\|f\|\|u\|$, and $\|f_\ii\|=\|f\|$. $\tr$ 

\medskip
As a corollary, for every ${\bf L}$--space $E$ and $u\in{\bf L}$ we have 
$\|u\|\ge\sup\{\|f_\ii(u)\|\}$, where supremum is taken over all $f\in E^*; \|f\|\le1$. But such a 
supremum is exactly $\|u\|_{\min}$. This justifies  the word ``minimal'' in Example 2.5.

\medskip 
Now we want to define amplifications of bilinear operators. For this we need a certain additional 
structure, connected with our base space, namely, a fixed bilinear operator, denoted in what 
follows by $\di:{\bf L}\times{\bf L}\to{\bf L}$ and called {\it $\di$--operation}. Let us write 
$\xi\di\eta$ instead of $\di(\xi,\eta)$. We call a $\di$--operation {\it metric}, if we always have 
$\|\xi\di\eta\|=\|\xi\|\|\eta\|$.  

If $E$ is a linear space, $\xi\in{\bf L}$ and $u\in{\bf L}E$, we set $\xi\di u:=T_\xi\cd u$, where 
$T_\xi\in\bb$ takes $\eta$ to $\xi\di\eta$. Thus, this version of a $\di$--operation is well 
defined on elementary tensors by the equality $\xi\di\eta  x:=(\xi\di\eta)x$. Similarly, we 
introduce $u\di\eta\in{\bf L}E$ by the equality $\xi x\di\eta:=(\xi\di\eta)x$. In the case of a 
metric $\di$--operation we obviously have $T_\xi=\|\xi\|S$, where $S$ is an isometry. Therefore, if 
$E$ is an ${\bf L}$--space, we have 
\[
\|\xi\di u\|=\|\xi\|\|u\|, \qq {\rm and\q similarly} \qq  \|u\di\eta\|=\|\eta\|\|u\|.\eqno(3.1)
\]

\medskip
{\bf Example 3.3}. Suppose that our base space is $L_p(X)$ from Example 2.2. We shall say that two 
given measure spaces are {\it of the same type}, if  a) they simultaneously have or do not have 
continuous ( = non-atomic) part, and b) they have sets of atoms of the same cardinality. As it is 
well known (see, e.g.,~\cite[Cor. 9.12.18]{bog}, and also~\cite[\S14]{lac} or~\cite[III.A]{woj}), 
in the case $p\ne2$ the spaces $L_p(X)$ and $L_p(Y)$ are isometrically isomorphic if, and only if 
the respective measure spaces are of the same type. 
% (а в противном случае они даже не топологически изоморфны). 
From this we immediately see that in the case $p\ne2$ the spaces $L_p(X)$ and $L_p(X\times X)$ are 
isometrically isomorphic if, and only if our $X$ is either devoid of atoms or has exactly one atom, 
or has infinite (necessarily countable) set of atoms. Thus, in these three cases the spaces 
$L_p(X\times X)$ and $L_p(X)$ are isometrically isomorphic. {\bf We choose an arbitrary
 isometric isomorphism $i:L_p(X\times X)\to L_p(X)$ and fix it throughout the whole 
paper}. After this we introduce $\di:L_p(X)\times L_p(X)\to L_p(X)$ as the composition 
$i\vartheta_{X}$, where the bilinear operator $\vartheta_{X}:L_p(X)\times L_p(X)\to L_p(X\times X)$ 
takes a pair $(x,y)$ to the function (of two variables) $x(s)y(t); s,t\in X$. Clearly, such a 
$\di$--operation is metric. 

\medskip
{\bf Remark 3.4}. Of course, if $p=2$, that is we deal with Hilbert spaces, a metric 
$\di$--operation exists for arbitrary (non-finite) measure spaces. This case was studied in details 
in~\cite{heh}. 

\medskip
{\bf From now on, and up to the end of the paper we assume that our base space is endowed with a 
metric $\di$--operation.} 

\medskip
%Вернемся к общим базовым пространствам, снабженным $\di$--операцией. 
Now let $\rr:E\times F\to G$ be a bilinear operator between linear spaces. Its {\it amplification} 
is the bilinear operator $\rr_\ii:{\bf L}E\times{\bf L}F\to{\bf L}G$, associated with the 4-linear 
operator ${\bf L}\times E\times{\bf L}\times F\to{\bf L}G:(\xi,x,\eta,y)\mt(\xi\di\eta)\rr(x,y)$. 
Thus, $\rr_\ii$ is well defined on elementary tensors by $\rr_\ii(\xi x,\eta 
y)=(\xi\di\eta)\rr(x,y)$ . 

\medskip
{\bf Definition 3.5.}. A bilinear operator $\rr$ between ${\bf L}$--spaces is called {\it ${\bf 
L}$--bounded}, respectively, {\it ${\bf L}$--contractive},  if its amplification is (just) bounded, 
respectively, contractive. We put $\|\rr\|_{{\bf L}b}:=\|\rr_\ii\|$. 

\medskip
Let $\rr$ be a ${\bf L}$--bounded bilinear operator. Then the equality $\rr_\ii(\xi x,\eta 
y)=(\xi\di\eta)\rr(x,y)$ implies that $\rr$,  being considered between respective underlying  
spaces is just bounded, and $\|\rr\|\le\|\rr\|_{{\bf L}b}$. At the same time, similarly to linear 
operators, sometimes the ``classical'' boundedness automatically implies the ${\bf 
L}$--boundedness. 

\medskip
{\bf Proposition 3.6.} {\it Let $E, F$ be ${\bf L}$-spaces, $f:E\to\co$ and $g:F\to\co$ bounded 
functionals. Then the bilinear functional $f\times g:E\times F\to\co:(x,y)\mt f(x)g(y)$ is ${\bf 
L}$--bounded, and $\|f\times g\|_{{\bf L}b}=\|f\|\|g\|$.}
 
\smallskip
$\tl$ Since $\|f\times g\|=\|f\|\|g\|$, it suffices to show that $\|f\times g\|_{{\bf 
L}b}\le\|f\|\|g\|$. Indeed, combining the obvious formula $(f\times g)_\ii(u,v)=f_\ii(u)\di 
g_\ii(v)$ with Proposition 3.2, we see that $\|(f\times g)_\ii(u,v)\|\le\|f\|\|g\|\|u\|\|v\|$.  
$\tr$ 

\medskip
In the following proposition $E$ is a normed space, $F$ is an ${\bf L}$--space. Let us denote the 
${\bf L}$--quantization of the space $E\ot_{pr}F$, considered in Example 2.4, by the same symbol
 $E\ot_{pr}F$; it will not lead to a misunderstanding. 

\medskip
{\bf Proposition 3.7.} {\it  The canonical bilinear operator $\vartheta:E_{\max}\times F\to E\ot_p 
F:(x,y)\mt x\ot y$, considered between the respective ${\bf L}$--spaces, is ${\bf L}$--contractive. 
} 

\smallskip
$\tl$ Consider the diagram 
\[
\xymatrix@C+20pt{{\bf L}E_{\max}\times{\bf L}F \ar[r]^{\vartheta_\ii}\ar[d]_{\beta_0\times\id_{{\bf L}F}}
& {\bf L}(E\ot_{pr} F) \ar[d]^{\beta} \\
(E\ot_{pr}{\bf L})\times{\bf L}F \ar[r]^{{\cal S}} & E\ot_{pr}{\bf L}F },
\]
\noindent where the bilinear operator ${\cal S}$ is well defined by taking $(x\ot\xi,v)$ to 
$x\ot(\xi\di v)$, $\beta$ is the isometric operator from Example 2.4, and the ``flip'' $\beta_0$ is 
its special case (when $F:=\co$). The diagram is obviously commutative, and ${\cal S}$ is 
contractive. As a corollary, 
%рутинные выкладки показывают, что 
for $w\in{\bf L}E_{\max}$ and $v\in{\bf L}F$ we have $\|\vartheta_\ii(w,v)\|\le\|w\|\|u\|$. $\tr$ 

\section{The general ${\bf L}$--tensor product}

Let us fix, throughout this section, two {\it arbitrary chosen} ${\bf L}$--spaces $E$ and $F$. 
Further, let $\mho$ be a subclass of the class of all normed ${\bf L}$-spaces. 

\medskip
{\bf Definition 4.1} A pair $(\Theta,\theta)$, consisting of $\Theta\in\mho$ and an ${\bf 
L}$--contractive bilinear operator $\theta:E\times F\to\Theta$, is called {\it tensor product of 
$E$ and $F$ relative to $\mho$} if, for every $G\in\mho$ and every ${\bf L}$--bounded bilinear 
operator $\rr:E\times F\to G$, there exists a unique ${\bf L}$--bounded operator $R:\Theta\to G$ 
such that the diagram 
\[
\xymatrix@R-10pt@C+15pt{
E\times F \ar[d]^{\theta} \ar[dr]^{\rr} & \\
\Theta \ar[r]^R  &  G  } %\eqno (4)
\]
\noindent is commutative, and moreover $\|R\|_{{\bf L}b}=\|\rr\|_{{\bf L}b}$. 

\medskip
Such a pair is unique in the following sense: if $(\Theta_k,\theta_k); k=1,2$ are two pairs, 
satisfying the given definition for a certain $\mho$, then there is a ${\bf L}$--isometric 
isomorphism $I:\Theta_1\to\Theta_2$, such that $I\theta_1=\theta_2$. This fact is a particular case 
of a general--categorical observation concerning the uniqueness of an initial object in a category; 
cf., e.g.,~\cite{mcl},~\cite[Theorem 2.73]{heb3}. However, the question about the existence of such 
a pair depends on our luck with the choice of the class $\mho$. 

\medskip
{\bf Definition 4.2}. The tensor product of $E$ and $F$ relative to the class of all normed ${\bf 
L}$--spaces is called {\it non-completed general ${\bf L}$--tensor product} of our spaces. 

\medskip
We shall prove the existence of such a pair, displaying its explicit construction. 

First, we need a sort of ``extended''  version of the diamond multiplication, this time between 
elements of amplifications of linear spaces. Namely, for $u\in{\bf L}E, v\in{\bf L}F$ we consider 
the element $u\di v:= \vartheta_\ii(u,v)\in{\bf L}(E\ot F)$, where $\vartheta:E\times F\to E\ot F$ 
is the canonical bilinear operator. In other words,  this ``diamond operation'' is well defined by 
$\xi x\di\eta y:=(\xi\di\eta)(x\ot y).$ 

%Let us make an obvious observation.

\medskip
{\bf Proposition 4.3.} {\it Every $U\in{\bf L}(E\ot F)$ can be represented as 
\[
U=\sum_{k=1}^na_k\cd(u_k\di v_k)\eqno(4.1)
\]
for some natural $n$ and $a_k\in\bb, u_k\in{\bf L}E,v_k\in{\bf L}F, k=1,...,n$. } 

\smallskip
$\tl$ Since every element of ${\bf L}(E\ot F)$ is the sum of several elements of the form $\xi(x\ot 
y);\xi\in{\bf L}, x\in E,y\in F$, it is sufficient to verify thr assertion on elements of the 
indicated form. Take an arbitrary vector $\eta\in{\bf L}$ and  an arbitrary operator $a\in\bb$, 
such that $a(\eta\di\eta)=\xi$. then we obviously have, $\xi(x\ot y)=a\cd(u\di v)$, where $u:=\eta 
x, v:=\eta y$. $ \tr$ 

\medskip
As a corollary, the operator $\bb\ot_{pr}{\bf L}E\ot_{pr}{\bf L}F\to{\bf L}(E\ot F):(a\ot u\ot 
v)\mt a\cd(u\di v)$ is surjective. Therefore the space ${\bf L}(E\ot F)$ can be endowed with the 
seminorm of the respective quotient space that we denote by $\|\cd\|_{{\bf L}}$. in other words, we 
have 
\[ 
\|U\|_{{\bf L}}:=\inf\{\sum_{k=1}^n\|a_k\|\|u_k\|\|v_k\|\}, \eqno(4.2)
\] 
where the infimum is taken over all possible representations of $U$ as indicated in (4.1). 

Being a quotient module of the module $\bb\ot_{pr}[{\bf L}E\ot_{pr}{\bf L}F]$, which is certainly 
contractive, the {\it seminormed $\bb$-module $({\bf L}(E\ot F),\|\cd\|_{{\bf L}})$ is itself 
contractive. } 
%Предложение 15. Полунормированный $\bb$--модуль $({\bf L}(E\ot F),\|\cd\|_{{\bf L}})$ -- сжимающий. 
Therefore $\|\cd\|_{{\bf L}}$ is an ${\bf L}$--seminorm on $E\ot F$. Denote the resulting
 ${\bf L}$--space by $E\ot_{\bf L}F$. 
%\smallskip
%$\tl$ Ясно, что $\bb\ot_p{\bf L}E\ot_p{\bf L}F$ -- это сжимающий левый $\bb$--модуль как тензорное 
%произведение левого $\bb$-модуля $\bb$ и линейного пространства ${\bf L}E\ot{\bf L}F$. Поэтому  
%${\bf L}(E\ot F)$ является образом сжимающего левого $\bb$-модуля при факторотображении 
%полунормированных пространств. Поскольку последнее -- это морфизм модулей, мы сразу получаем 
%желаемое свойство. $\tr$ 

Observe the obvious estimate
\[
\|u\di v\|_{{\bf L}}\le\|u\|\|v\|; \q u\in{\bf L}E,v\in{\bf L}F.\eqno (4.3) 
\]

Since $u\di v=\vartheta_\ii(u,v)$, we see that {\it $\vartheta$, being considered with values in 
$E\ot_{{\bf L}}F$, is ${\bf L}$--contractive.} 

Using that for $\xi\in\bf L;\|\xi\|=1$ we have $\|x\ot y\|=\|(\xi\di\xi)x\ot y\|$, we obtain in the 
underlying space of $E\ot_{\bf L}F$ the estimate 

\[
\|x\ot y\|\le\|x\|\|y\|; \q x\in E, y\in F.\eqno (4.4)
\]

%(На самом деле в (4.3) и (4.4) имеют место равенства, но мы не будем здесь это обсуждать). 
\medskip
{\bf Proposition 4.4}. {\it Let $G$ be a ${\bf L}$-space, $\rr:E\times F\to G$ an ${\bf 
L}$--bounded bilinear operator, $R:E\ot_{\bf L}F\to G$ the associated linear operator. Then $R$ is 
${\bf L}$--bounded, and $\|\rr\|_{{\bf L}b}=\|R\|_{{\bf L}b}$. } 

\smallskip
$\tl$ Take $U\in{\bf L}(E\ot_{{\bf L}}F))$ and represent it according to (4.1). Since $R_\ii$ is a 
$\bb$-module morphism, we have, using the obvious equality $R_\ii(u\di v)=\rr_\ii(u,v)$, that 
$R_\ii(U)=\sum_{k=1}^na_k\cd\rr_\ii(u_k, v_k)$. Consequently, we have 
\[ 
\|R_\ii(U)\|\le\sum_{k=1}^n\|a_k\|\|\rr_\ii(u_k, 
v_k)\|\le\|\rr\|_{{\bf L}b}\sum_{k=1}^n\|a_k\|\|u_k\|\|v_k\|. 
\]
From this, using (4.2), we  obtain that $\|R\|_{{\bf L}b}\le\|\rr\|_{{\bf L}b}$. The inverse 
inequality follows from (4.3). $\tr$ 

\medskip
{\bf Proposition 4.5}. {\it {\rm (As a matter of fact),} $\|\cd\|_{{\bf L}}$ is a norm. }

\smallskip
$\tl$  By Proposition 2.4, it is sufficient to show that, for a non-zero elementary tensor $\xi w; 
w\in E\ot_{{\bf L}}F, \xi\in{{\bf L}}; \|\xi\|=1,$ we have $\|\xi w\|_{{\bf L}}\ne0$. Since 
  $E$ and $F$ are normed spaces, then there exist bounded functionals $f:E\to\co, g:F\to\co$ such that 
  $(f\ot g)w\ne0$. Set in Proposition 4.4 $\rr:=f\times g:E\times F\to\co$. By virtue of Proposition 3.6,
   $\rr$ is ${\bf L}$--bounded, hence the operator $(f\ot g)_\ii:{\bf L}(E\ot_{\bf L}F)\to{\bf L}$ is 
bounded. At the same time $(f\ot g)_\ii(\xi w)=[(f\ot g)(w)]\xi\ne0$. $\tr$

\medskip
Combining the accumulated facts, we immediately obtain the desired existence theorem:
 
\medskip
{\bf Theorem 4.6.} {\it The pair $(E\ot_{\bf L}F,\vartheta)$ is a non-completed general ${\bf 
L}$--tensor product of $E$ and $F$. } 

\medskip
For some concrete tensor factors the introduced tensor product also becomes something concrete: 
%and ``palpable'': 
 
\medskip
{\bf Theorem 4.7}. {\it Let $E$ and $F$ be the spaces from Example 2.6. Suppose that ${\bf 
L}:=L_p(X)$, where $X$ satisfies the conditions, indicated in Example 3.3, and the $\di$-operation 
is taken from the same example. Then there exists an ${\bf L}$--isometric isomorphism
 $I:E_{\max}\ot_{{\bf L}}F\to E\ot_{pr}F$, acting as the identity operator on the common underlying
 linear space of our ${\bf L}$--spaces.} 

\smallskip
$\tl$ Consider the operator $I$, associated with $\vartheta$ (cf. Proposition 3.7). It acts as it 
is indicated in the formulation and, by virtue of Theorem 4.6, it is ${\bf L}$--contracting. We 
must only show that $I_\ii$ does not decrease norms of vectors. 

Take $U\in{\bf L}(E\ot F)$. Identifying ${\bf L}(E\ot F)$ with $E\ot{\bf L}F$, we can represent $U$ 
as $\sum_{k=1}^n x_k\ot v_k; x_k\in E, v_k\in{\bf L}F$. 

Let $p<\ii$ and $q$ its conjugate. Then we choose an arbitrary $e\in{\bf L};\|e\|=1$, denote by 
$e^*\in L_q(X)$ a function of norm 1, for which we have $\int_Xe(s)e^*(s)ds=1$, and consider the 
operator $j:L_p(X\times X)\to L_p(X)$, taking a function $f(s,t)$ to $g(t):=\int_Xf(s,t)e^*(s)ds$. 
If $p=\ii$, we set $e(s)\equiv1$ and introduce $j:L_p(X\times X)\to L_p(X)$, taking  $f(s,t)$ to 
$g(t):=ess\;sup|f^t|\textbf{}$, where $f^t(s):=f(s,t)$. Then in both cases we set 
$T:=ji^{-1}\in\bb$ and see that $\|T\|=1$. Moreover, representing every $v_k$ as the sum of 
elementary tensors from ${\bf L}F$, we easily obtain that 

$$
U=T\cd\left[\sum_{k=1}^nex_k\di v_k \right].
$$
 From this, by virtue of (4.2), we obtain the 
estimate $\|U\|_{{\bf L}b}\le\sum_{k=1}^n\|x_k\|\|v_k\|$ and, as a corollary, we have 
\[
\|U\|_{\bf L}\le\inf\{\sum_{k=1}^n\|x_k\|\|v_k\|\},\eqno (4.5)
\] 
where the infimum is taken over all representations of $U$ in the indicated form.

Now look at $I_\ii(U)$. It is the same $\sum_{k=1}^n x_k\ot v_k$, only considered in the normed 
space $E\ot_{pr}{\bf L}F$. It follows that $\|I_\ii(U)\|$ is exactly the infimum, indicated in 
(4.5). Thus, the desired estimate  $\|I_\ii(U)\|\ge\|U\|_{{\bf L}}$ is obtained.  $\tr$ 

\medskip
{\bf Remark 4.8.} As an easy corollary of this theorem, we have, up to an ${\bf L}$--isometric 
isomorphism, that $E_{\max}\ot_{{\bf L}} F_{\max}=[E\ot_{pr} F]_{\max}$.
 In particular, for a Hilbert space $H$ we have 
$H_{\max}\ot_{{\bf L}}H_{\max}={\cal N}_0(H)_{\max}$, where ${\cal N}_0(H)$ is the space of finite 
rank operators on $H$, equipped with the trace class norm. 
 
\section{$p$--convexity and  $p$--convex tensor product}

Now we need one more, apart from $\di$--operation, additional structure in our base space. Namely, 
we say that ${\bf L}$ is {\it a stratified space}, if a certain family ${\cal P}$ of projections of 
norm 1 (or 0), acting on ${\bf L}$, is distinguished, and it is such that $P,Q\in{\cal P}$ implies 
$PQ=QP\in{\cal P}$, and if $P,Q\in{\cal P}$ are orthogonal, then $P+Q\in{\cal P}$. Projections that 
belong to ${\cal P}$ will be called {\it proper.} 
%shall speak about {\it a space with projections}, if in ${\bf L}$ a certain family ${\cal P}$ of 
%projections of norm 1 (or 0) is distinguished such $P,Q\in{\cal P}$ implies $PQ=QP\in{\cal P}$, and 
%if $P,Q\in{\cal P}$ are orthogonal, then $P+Q\in{\cal P}$. Projections that belong to ${\cal P}$ 
%will be called {\it proper.} 

\medskip
{\bf Example 5.1.} Let ${\bf L}:=L_p(X)$. If $X'$ is a measurable subset in $X$, we denote by 
$P_{X'}\in\bb(L_p(X))$ the projection, acting as $f\mt f\chi$, where $\chi$ is the characteristic 
function of $X'$. Of course, if the measure of $X'$ is positive, we have $\|P_{X'}\|=1$. We shall 
identify the image of this projection with $L_p(X')$. Clearly, projections of this sort are 
orthogonal if, and only if the respective sets have intersection of measure 0. It is obvious that 
the family of projections of the indicated form satisfies the conditions, formulated above. 
Speaking about $L_p(X)$ as of a stratified space, we shall always mean this particular family. 

\medskip
{\it In what follows, if numbers $\lm_k\ge0; k=1,...,n$ are given, we shall understand the 
expression $(\sum_{k=1}^n\lm_k^p)^{\frac{1}{p}}$ as $\max\{\lm_1,...,\lm_n\}$ in the case $p=\ii$.}

Let us distinguish, for convenience of future references, the following obvious observation.

\medskip
{\bf Proposition 5.2.} {\it Let $X$ and $Y$ be two measure spaces, $S_k\in\bb(L_p(X), L_p(Y)); 
\\ k=1,\dots,n$, $X_k'$ и $Y_k'$ -- two families of pairwise non-intersect measurable subsets in  $X$ 
и $Y$, respectively. Then we have 
$\|\sum_{k=1}^nP_{Y_k'}S_kP_{X_k'}\|\le\max\{\|S_k\|;k=1,\dots,n\}$.} $\tl\q\tr$ 

\medskip
Let $u$ be an element of a certain seminormed ${\bf L}$--space. We call a projection $P\in\bb$ a 
{\it support } of $u$, if $P\cd u=u$. 

\medskip
{\bf Definition 5.3.} Let ${\bf L}$ be a stratified space, $1\le p\le\ii$. Тhen a seminormed ${\bf 
L}$--space $E$ is called {\it $p$--convex,} if for every $u,v\in{\bf L}E$ with orthogonal proper 
supports, we have $\|u+v\|\le(\|u\|^p+\|v\|^p)^{\frac{1}{p}}$. As an immediate corollary, for 
$u_1,...,u_n\in{\bf L}E$ with pairwise orthogonal supports from ${\cal P}$, we have 
$\|\sum_{k=1}^nu_k\|\le(\|\sum_{k=1}^n\|u_k\|^p)^{\frac{1}{p}}$. 

\medskip
For the special case ${\bf L}:=\ell_p$ this definition is equivalent to the definion of a 
$p$--convex $p$--multi-normed space in~\cite{dal}. Also it worthy to mention, in this connection, 
the theory of $p$--operator spaces of Daws~\cite{daw}; see also earlier papers of Pisier~\cite{pi2} 
and Le Merdy~\cite{lem}. 

As to the base space ${\bf L}$ itself, it is called $p$--convex, if it is $p$--convex after its 
identification with the ${\bf L}$--space $\co$. Needless to say, $L_p(X)$, as a base space, is
$p$--convex.

Clearly, for every stratified ${\bf L}$ all seminormed ${\bf L}$--spaces are 1-convex. Also it is 
obvious that every $p$--convex space is $r$--convex for each $r; 1\le r<p$. 

\medskip 
{\bf Proposition 5.4.} {\it If ${\bf L}$ is $p$--convex, in particular, if ${\bf L}:=L_p(X)$, then 
every ${\bf L}$--space $E$ with the minimal quantization is  $р$-convex.} 

\smallskip
$\tl$ If $u,v\in{\bf L}E$ have orthogonal supports, then for every $f\in E^*$ elements 
$f_\ii(u),f_\ii(v)$ have orthogonal supports in ${\bf L}$. Therefore we have 
$\|f_\ii(u+v)\|^p\le\|f_\ii(u)\|^p+\|f_\ii(v)\|^p$. It remains to take the relevant supremum in the 
right part of the inequality over all $f\in E^*; \|f\|=1$, and then do the same with the left part. 
$\tr$ 
 
\medskip
At the same time in the case $p>1$ and ${\bf L}:=L_p(X)$ the maximal quantization of a normed space 
is not, generally speaking, $p$--convex. The space $E:=\ell_1$ serves as the simplest 
counter-example. 

One can construct $p$--convex spaces, embarking from other $p$--convex spaces. For example, it is 
not difficult to show that, for $q_k$--convex ${\bf L}$--spaces $E_k;k=1,2$ and 
$p\le\min\{q_1,q_2\}$, the ${\bf L}$--space $E_1\oplus E_2$, being considered with the ${\bf 
L}$--norm of the $\ell_p$--sum of normed spaces ${\bf L}E_1$ and ${\bf L}E_2$, is $p$--convex. 

%Вот пример. 
%Пусть $({\bf L},{\cal P})$ -- пространство c проекторами, $E_k; k=1,2$ -- два ${\bf 
%L}$--пространства. Рассмотрим ${\bf L}$--пространство $E_1\oplus E_2$ с ${\bf L}$--нормой 
%$\ell_p$--суммы нормированных пространств ${\bf L}E_1$ и ${\bf L}E_2$ (Такое пространство является 
%${\bf L}$--квантованием пространства $E_1\oplus_p E_2$). 

%\medskip
%Предложение 1. Пусть ${\bf L}$--пространство $E_k; k=1,2$ $q_k$--выпукло, где $1\le q_1,q_2\le\ii$. 
%Тогда при $p\le\min\{q_1,q_2\}$ ${\bf L}$--пространство $E_1\oplus_p E_2$ $p$--выпукло. 

%\smallskip
%$\tl$  Пусть $u,v\in E_1\oplus E_2$ обладают ортогональными носителями $P,Q$. Очевидно, $u=u_1+u_2, 
%v=v_1+v_2$, где $u_k,v_k\in E_k; k=1,2$. Тогда $P$ -- носитель $u_1$ и $u_2$, а $Q$ -- носитель   
%$v_1$ и $v_2$. Далее, $\|u+v\|=(\|u_1+v_1\|^p+\|u_2+v_2\|^p)^{\frac{1}{p}}$. Но 
%$\|u_k+v_k\|^{q_k}\le\|u_k\|^{q_k}+\|v_k\|^{q_k}$, а это влечет 
%$\|u_k+v_k\|^{p}\le\|u_k\|^{p}+\|v_k\|^{p}$. Следовательно, 
%$\|u+v\|\le(\|u_1\|^p+\|v_1\|^p+\|u_2\|^p+\|v_2\|^p)^{\frac{1}{p}}=(\|u\|^p+\|v\|^p)^{\frac{1}{p}}$. $\tr$ 

\medskip
Now suppose that we are given a space ${\bf L}$, which is stratified and has a $\di$--operation. 
Fix two ${\bf L}$--spaces $E$ and $F$. 

\medskip
{\bf Definition 5.5.} The tensor product of $E$ and $F$ relative to the class of all $p$--convex 
${\bf L}$--spaces is called {\it non-completed $p$--convex ${\bf L}$--tensor product} of our 
spaces. 

\medskip
Unfortunately, at the moment we can prove the existence of such a tensor product only under rather 
burdensome additional assumptions on $\di$ and especially on ${\cal P}$. In fact, all our
%now we do not have 
examples of triples $({\bf L},{\cal P},\di)$, satisfying these conditions, are, in a sense, 
%that would not be 
too close to triples, arising in the case of base spaces  $L_p(X)$, moreover with (however mild) 
assumptions on $X$. Therefore we decided not to bother the reader with the list of these 
conditions. Instead, 
% Так что, не мудрствуя лукаво, 

{\bf from now on up to the end of the paper we consider, in capacity of base spaces, only spaces 
$L_p(X)$, where $X$ is a measure space that either has no atoms or has infinite set of atoms.} 

We call a measure space of the indicated sort {\it convenient}. The family of proper projections, 
acting on $L_p(X)$, is defined according to Example 5.1, and a $\di$--operation according to 
Example 3.3. Note the the case of a single atom, permitted in the latter example, now is forbidden; 
otherwise (here we open our cards) the future Proposition 5.6 fails to be true. 

We call an isometric operator on $L_p(X)$ {\it proper}, if its image coincides with the image of 
some proper projection or, equvalently, its image is $L_p(X')$ for some measurable subset $X'$ of $ 
X$. We call two proper isometries {\it disjoint}, if the intersection of their images is $\{0\}$ 
or, equvalently, the corresponding projections are orthogonal. Since $X$ is convenient, it contains 
an infinite family of pairwise disjoint measurable subsets of the same type as $X$ (cf. Example 
3.3). Here is an immediate corollary. 

\medskip
{\bf Proposition 5.6. } {\it There exists an infinite family of pairwise disjoint proper 
isometries, acting on $L_p(X)$.} $\tl \q \tr$ 

\medskip
(By the way, in the case of general ${\bf L}$ it would be sufficient to find two disjoint proper 
isometries, say $I_1$ and $I_2$. Then the isometries $I_1,I_1I_2,I_1I_2^2,I_1I_2^3,...$ would fit). 

\medskip
If $I\in\bb$ is a proper isometry with image $L_p(X')$, we denote by $I^\star\in\bb$ the operator 
$I^{-1}P_{X'}$, that is the coisometric operator ( = quotient map), acting as $I^{-1}$ on $L_p(X')$ 
and vanishing on the complementary subspace $L_p(X\setminus X')$. It is clear that, if $I_k; 
k=1,...,n$ are pairwise disjoint proper isometries, we have 
\[
I_k^\star I_l=\de^k_l\id_{\bf L}. \eqno (5.1)
\]

We proceed to the explicit construction of the $p$-convex tensor product. 

\medskip
{\bf Proposition 5.7}. {\it If $E,F$ are linear spaces, then every $U\in L(E\ot F)$ can be 
represented as 
\[
a\cd\sum_{k=1}^n I_k\cd(u_k\di v_k),\eqno (5.2)
\]
where $a\in\bb$, and $I_k$ are pairwise disjoint proper isometries. } 
 
\smallskip
$\tl$ Represent $U$ as in (4.1). By virtue of Proposition 5.6, there exist $n$ proper pairwise 
disjoint isometries $I_k$. Consider in ${\bf L}(E\ot F)$ the element 
  \[
 \left(\sum_{k=1}^na_kI^\star_k\right)\cd\left(\sum_{l=1}^nI_l\cd(u_l\di v_l)\right).
 \]
It follows from (5.1) that it is exactly $U$. $\tr$ 

\medskip
From now on we assume that we are given two arbitrary (not necessarily $p$--convex!) ${\bf 
L}$--spaces $E$ and $F$. Assign to every $U\in L(E\ot F)$  the number 
 \[
\|U\|_{p{\bf L}}:=\inf\left\{\|a\|\left(\sum_{k=1}^n\|u\|^p\|v\|^p\right)^{\frac{1}{p}}\right\},\eqno (5.3)
\]
where the infimum is taken over all possible representations of $U$ in the form (5.2). 
%(From now on we suppose that $p<\ii$; the case $p=\ii$, being in a hurry, in this draft we leave to the %reader).

We distinguish the obvious 

\medskip
{\bf Proposition 5.8}. {\it For every $U\in{\bf L}(E\ot F)$ and $a\in\bb$ we have $\|a\cd 
U\|_{p{\bf L}}\le\|a\|\|U\|_{p{\bf L}}$. } $\tl \q \tr$ 

\medskip
What is less obvious, it is 
 
\medskip
{\bf Proposition 5.9}. {\it The function $U\mt\|U\|_{{p{\bf L}}}$ is a seminorm on ${{\bf L}}(E\ot 
F)$. } 

\smallskip
$\tl$ Suppose that $U=a\cd\sum_{k=1}^nI'_{k}\cd(u^1_k\di v^1_k)$ and $V=b\cd\sum_{l=1}^m 
I''_l\cd(u^2_l\di v^2_l)$, where $I'_{k}$ are proper pairwise disjoint isometries, and the same is 
true for $I''_{l}$. 

Choose, in an arbitrary way, one more pair of disjoint isometries $I_U,I_V\in\bb$ and observe that 
$$
U+V=\left(aI_U^\star+bU_V^\star\right)\cd
\left(\sum_{k=1}^n I_UI'_{k}\cd(u^1_k\di v^1_k)+\sum_{l=1}^mI_VI''_l\cd(u^2_l\di v^2_l)\right).
$$
Evidently, the compositions $I_UI'_k$ and $I_VI''_l$ are proper isometries and, being considered 
all together, they are pairwise disjoint. Therefore, by virtue of (5.3) and the previous 
proposition, we have 
$$
\|U+V\|_{p{\bf L}}\le\|aI_U^\star+bI_V^\star\|
\left(\sum_{k=1}^n\|u^1_k\|^p\|v^1_k\|^p+\sum_{l=1}^m\|u^2_l\|^p\|v^2_l\|^p\right)^{\frac{1}{p}}.\eqno(5.4)
$$

Our nearest aim is to obtain the estimate
$\|aI_U^\star+bI_V^\star\|\le(\|a\|^{q}+\|b\|^{q})^{\frac{1}{q}}$, where $q$ 
%:=p/(p-1)$ 
is the conjugate number to $p$. (As usual, we hold that 1 and $\ii$ are mutually conjugate). 

Take $\xi\in{\bf L},\|\xi\|\le1$ and denote the proper projections, corresponding to our 
isometries, by $P_U=I_UI_U^\star$ and $P_V=I_VI_V^\star$. Then, taking into account the G\"older 
inequality, we obtain that 
\[
\|(aI_U^\star+bI_V^\star)(\xi)\|\le\|(aI_U^\star P_U(\xi)\|+\|(bI^\star_VP_V)(\xi)\|\le
%\|(aS^{-1}P_H)(\xi)\|+\|(bT^{-1}P_K)(\xi)\|\le
\]
\[
%\|a\|\|\|I^\star_U\|\|P_U(\xi)\|+\|b\|\|\|I^\star_V\|\|P_V(\xi)\|\|\le
\|a\|\|\|P_U(\xi)\|+\|b\|\|\|P_V(\xi)\|\le
(\|a\|^{q}+\|b\|^{q})^{\frac{1}{q}}(\|P_U(\xi)\|^p+\|P_V(\xi)\|^p)^{\frac{1}{p}}
\]
%Ввиду неравенства Гёльдера, последнее число 
%\[
%\le(\|a\|^{q}+\|b\|^{q})^{\frac{1}{q}}(\|P_U(\xi)\|^p+\|P_V(\xi)\|^p)^{\frac{1}{p}},
%\]
Since our projections are orthogonal, the scond factor is equal to 
$((\|P_U(\xi)+P_V(\xi)\|)^{p})^{\frac{1}{p}}$, and therefore it does not exceed $\|P_U+P_V\|=1$. 
%\[
%\|P_U(\xi)+P_V(\xi)\|^{p\frac{1}{p}}\le\|P_u+P_V\|=1.
%\]
Our desired estimate follows. 

\medskip
Obviously, we can obtain representations of our $U$, by multiplying $a$ on a certain constant and 
dividing all $u_k^1$ on the same constant; in a similar way we can deal with $V$. Consequently, we 
have a right to assume in the case $1<p<\ii$ that $\|a\|^{q}=\sum_{k=1}^n\|u^1_k\|^p\|v^1_k\|^p$ 
and $\|b\|^{q}=\sum_{l=1}^m\|u^2_l\|^p\|v^2_l\|^p$. Also we can assume in the case $p=1$ that 
$\|a\|=\|b\|=1$, and in the case $p=\ii$ that $\max\{\|u^1_k\|\|v^1_k\|; 
k=1,...,n\}=\max\{\|u^2_l\|\|v^2_l\|; l=1,...,m\}$. 

Consequently, in the case $1<p<\ii$ we have 
\[
\|U+V\|_{p{\bf L}}\le(\|a\|^{q}+\|b\|^{q})^{\frac{1}{q}}(\|a\|^{q}+\|b\|^{q})^{\frac{1}{p}}=
\|a\|^{q}+\|b\|^{q}.
\]
But $q-1=q/p$, hence 
 $\|a\|^{q}=\|a\|(\sum_{k=1}^n\|u^1_k\|^p\|v^1_k\|^p)^{\frac{1}{p}}$. Together with the similar equality
 for $\|b\|^{q}$, this provides the estimate
 \[
\|U+V\|_{{p{\bf L}}}\le(\|a\|\left(\sum_{k=1}^n\|u^1_k\|^p\|v^1_k\|^p\right)^{\frac{1}{p}}+
\|b\|\left(\sum_l^m\|u^2_l\|^p\|v^2_l\|^p\right)^{\frac{1}{p}}. \eqno(5.5)
\]
It is easy to verify that the same estimate is valid in the remaining cases. 
%Если $p=1$, то $\|U+V\|_{p{\bf 
%L}}\le\max\{\|u\|,\|v\|\}\left(\sum_{k=1}^n\|u^1_k\|\|v^1_k\|+\sum_{l=1}^m\|u^2_l\|\|v^2_l\|\right)$, 
%а если $p=\ii$, то $\|U+V\|_{p{\bf L}}\le(\|u\|+\|v\|)\max\{\|u^1_k\|\|v^1_k\|,\|u^2_l\|\|v^2_l\|; 
%k=1,...,n,l=1,...,m\}$. Поэтому снова выполнена оценка (5.5). 
 
Thus, in all cases, taking the infimum from (5.3), we obtain the triangle inequality: 
$\|U+V\|_{p{\bf L}}\le\|U\|_{p{\bf L}}+\|V\|_{p{\bf L}}$. 

The property of seminorms, concerning the scalar multiplication, is immediate. $\tr$ 

\medskip
We see that $\|\cd\|_{\bf L}$ is an ${\bf L}$--seminorm on $E\ot F$. Denote the resulting ${\bf 
L}$--seminormed space by $E\ot_{p{{\bf L}}}F$. 

\medskip
{\bf Proposition 5.10}. {\it The introduced space is $p$--convex. } 

\smallskip+
$\tl$ Let $U,V\in{\bf L}(E\ot F)$ have orthogonal supports $P_1$ and $P_2$. Choose their arbitrary 
representations in the same form, as in Proposition 5.9, and take corresponding $I_U$ and $I_V$. 
The estimate (5.4) appears. Obviously, we have a right to assume that $a=P_1a,b=P_2b$, and also 
that  $\|a\|=\|b\|=1$. Thus, in the notations $P_U:=I_UI_U^\star$ and $P_V:=I_VI_V^\star$ we have 
$aI'^\star=P_1aI'^\star P_U$ и $bI''^\star=P_2bI''^\star P_V$. From this, by virtue of Proposition 
5.2, we have 
$\|aI'^\star+bI''^\star\|\le\max\{\|aI'^\star\|,\|bI''^\star\|\}\le\max\{\|a\|,\|b\|\}=1$, hence  
$\|U+V\|_{p{\bf L}}\le\left(\sum_{k=1}^n\|u^1_k\|^p\|v^1_k\|^p+ 
\sum_{l=1}^m\|u^2_l\|^p\|v^2_l\|^p\right)^{\frac{1}{p}}$. It remains to take the corresponding 
infima. $\tr$ 

\medskip
Similarly to the case of the tensor product `` $\ot_{\bf L}$ '', we have the estimate 
\[
\|u\di v\|_{p{\bf L}}\le\|u\|\|v\|; \q u\in{\bf L}E,v\in{\bf L}F.  \eqno(5.6)
\]
It follows that {\it the canonical bilinear operator $\vartheta:E\times F\to E\ot_{p{{\bf L}}}F$} 
%(отличающийся от рассмотренного в Разделе 4 тем, что он принимает значения в другом 
%нормированном пространстве) 
(just as $\vartheta:E\times F\to E\ot_{{\bf L}}F$ in Section 4) {\it is ${{\bf L}}$--contractive.} 

The same argument that provides the estimate (4.4), shows that in the underlying seminormed space 
of $E\ot_{p{\bf L}}F$ we have the estimate 
\[
\|x\ot y\|\le\|x\|\|y\|; \q x\in E, y\in F.\eqno (5.7).
\]

\medskip
{\bf Proposition 5.11}. {\it Let $G$ be a $p$--convex ${\bf L}$--space, $\rr:E\times F\to G$ an 
${\bf L}$--bounded bilinear operator. Then the associated linear operator $R:E\ot_{p{\bf L}}F\to G$ 
is ${\bf L}$--bounded, and   $\|\rr\|_{{\bf L}b}\|=\|R\|_{{\bf L}b}$. } 

\smallskip
Take $U\in{\bf L}(E\ot_{l}F))$ and its representation as in (5.2). Since $R_\ii$ is a $\bb$--module 
morphism, we see that $R_\ii(U)=a\cd(\sum_{k=1}^nI_k\cd\rr_\ii(u_k, v_k))$. Look at the elements 
$I_k\cd\rr_\ii(u_k,v_k)\in{\bf L}G$. They have pairwise orthogonal supports, namely $I_kI_k^\star$, 
and $G$ is $p$--convex. From this we obtain that 
$$
\|R_\ii(U)\|\le\|a\|\left(\sum_{k=1}^n\|I_k\cd\rr_\ii(u_k,v_k)\|^p\right)^{\frac{1}{p}}\le
\|a\|\left(\sum_{k=1}^n\|\rr_\ii(u_k,v_k)\|^p\right)^{\frac{1}{p}}\le
$$
$$
\|a\|\left(\sum_{k=1}^n\|\rr\|_{b{\bf L}}^p\|u_k\|^p\|v_k\|^p\right)^\frac{1}{p}=
\|\rr\|_{b{\bf L}}\|a\|\left(\sum_{k=1}^n\|u_k\|^p\|v_k\|^p\right)^\frac{1}{p}.
$$
Therefore $\|R\|_{{b{\bf L}}}\le\|\rr\|_{{b{\bf L}}}$. The inverse inequality follows from (5.6).  
$\tr$ 

\medskip
{\bf Proposition 5.12}. {\it {\rm (As a matter of fact),} $\|\cd\|_{p{\bf L}}$ is a norm. } 

\smallskip
$\tl$ Needless to say that $\co$ is a $p$--convex ${\bf L}$--space, hence we have a right to use 
Proposition 5.11. Therefore the proof of Proposition 4.5 goes up to obvious modifications. $\tr$

\medskip 
Combining the relevant propositions, we obtain the following existence theorem: 
 
\medskip
{\bf Theorem 5.13.} {\it The pair $(E\ot_{p{\bf L}}F,\vartheta)$ is a non-completed $p$--convex 
tensor product of ${\bf L}$--spaces $E$ and $F$. } 

\medskip
Such a theorem, in our opinion, can be considered as a far-reaching generalization of certain 
results of Lambert concerning the maximal tensor product of his ``Operatorfolgenr\"aume''; 
see~\cite[pp. 73-78]{lam}. In this connection we recall the papers of Blecher/Paulsen~\cite{blp} 
and Effros/Ruan~\cite{er2} about projective tensor products of operator spaces: they are at the 
source of all these constructions. 

\medskip
{\bf Remark 5.14}. We do not discuss here the ``non-discrete'' version of another, the so-called 
minimal tensor product, that was introduced by Lambert (in the frame-work of the ``coordinate'' 
approach) for 2--convex $l_2$--spaces in~\cite[3.1.3]{lam}. 

\section{$p$--convex tensor product of spaces $L_q(\cd)$}

In conclusion, we shall present an example, 
%show the situation, 
when the  $p$--convex tensor product 
%, that was introduced above for ${\bf L}:=
%$L_p(\cd)$--spaces 
%and {\it convenient} $X$, 
acquires especially transparent form.   It happens that in the case $1<p<\ii$ one should take, in 
the role of ``the best'' tensor factors, 
% in a sense, tensor factors are 
the spaces $L_q(\cd); q:=p/(p-1)$ with the minimal quantization, discussed in Example 2.5. 

\medskip
We remember that our base space ${\bf L}$ is $L_p(X)$ for convenient $X$. Moreover, {\bf throughout 
this section, all ${\bf L}$--spaces are supposed to be endowed with the minimal quantization}. 

\medskip
Let $Y$ and $Z$ be two measure spaces. Consider the linear operator \\ $J:L_p(Y)\ot L_p(Z)\to 
L_p(Y\times Z)$, well defined by the rule $x\ot y\mt x(s)y(t); \\ s\in Y, t\in Z$. It is easy to 
see that it is injective, and its image consists of functions of the form  
$\sum_{k=1}^nf_k(s)g_k(t); f_k\in L_p(Y),g_k\in L_p(Z)$. We see that this image is a normed 
subspace in $L_p(Y\times Z)$, which is dense provided $p<\ii$. Denote it by  $L_p(Y)\ot^p L_p(Z)$. 
Obviously, we can identify this space with the tensor product $L_p(Y)\ot L_p(Z)$, endowed with the 
respective induced norm. 

\medskip
{\bf Proposition 6.1.} {\it Let $A:L_p(Y_1)\to L_p(Y_2)$ and $B:L_p(Z_1)\to L_p(Z_2)$ be two 
bounded operators. Then the operator $A\ot B:L_p(Y_1)\ot^p L_p(Z_1)\to L_p(Y_2)\ot^p L_p(Z_2)$ is 
also bounded, and $\|A\ot B\|\le\|A\|\|B\|$.} 

\smallskip
$\tl$ Every $u\in L_p(Y_1)\ot^p L_p(Z_1)$ is a function of the form 
 $\sum_{k=1}^nf_k(s)g_k(t): \\ f_k\in L_p(Y_1), g_k\in L_p(Z_1)$. If $p<\ii$, then, by virtue of 
 Fubini Theorem we have
 $$
\|A\ot B(u)\|=
%\left(\int_{Y_2}\int_{Z_2}\left|\sum_kAf_k(s)Bg_k(t)\right|^pd(s,t)\right)^{\frac{1}{p}}=
\left(\int_{Z_2}\left[\int_{Y_2}\left|A(\sum_{k=1}^n((Bg_k)(t)f_k)(s))\right|^pds\right]dt\right)^{\frac{1}{p}}\le
$$
$$
%\left(\int_{Z_2}\left\|A(\sum_k((Bg_k)(t)f_k))\right\|_{L_p(Y_2)}^pdt\right)^{\frac{1}{p}}\le
\|A\|\left(\int_{Z_2}\left\|\sum_{k=1}^n(Bg_k)(t)f_k\right\|_{L_p(Y_2)}^pdt\right)^{\frac{1}{p}}=
\|A\|\left(\int_{Z_2}\left[\int_{Y_2}\left|\sum_{k=1}^n(Bg_k)(t)f_k(s)\right|^pds\right]dt\right)^{\frac{1}{p}}=
$$
%$$
%\|A\|\left(\int_{Z_2}\left[\int_{Y_2}\left|\sum_k(Bg_k)(t)f_k(s)\right|^pds\right]dt\right)^{\frac{1}{p}}=
%\|A\|\left(\int_{Y_2}\left[\int_{Z_2}\left|\sum_k(Bf_k(s)g_k)(t)\right|^pdt\right]ds\right)^{\frac{1}{p}}=
%$$
$$
\|A\|\left(\int_{Y_2}\left\|\sum_{k=1}^n(Bf_k(s)g_k)\right\|^p_{L_p(Z_2)}ds\right)^{\frac{1}{p}}\le
%\|A\|\left(\int_{Y_2}\left\|B(\sum_k(f_k(s)g_k)\right\|_{L_p(Z_2)}^pds\right)^{\frac{1}{p}}\le
\|A\|\|B\|\left(\int_{Y_2}\left\|\sum_{k=1}^n(f_k(s)g_k)\right\|_{L_p(Z_2)}^pds\right)^{\frac{1}{p}}=
$$
%$$
%\|A\|\|B\|\left(\int_{Y_2}\left\|\sum_k(f_k(s)g_k)\right\|_{L_p(Z_2)}^pds\right)^{\frac{1}{p}}=
%\|A\|\|B\|\left(\int_{Y_2}\left\|\sum_k(f_k(s)g_k(t))\right\|^pds\right)^{\frac{1}{p}}=
%$$
$$
\|A\|\|B\|\left(\int_{{Y_2}\times Z_2}\left|\sum_{k=1}^n(f_k(s)g_k(t))\right|^pd(s,t)\right)^{\frac{1}{p}}=
\|A\|\|B\|\|u\|. 
$$
If $p=\ii$, then a similar argument works, 
%To prove the assertion in the case  we use a similar argument, 
only instead of Fubini Theorem we apply, for functions $h\in L_\ii(Y,Z)$, $h_s\in L_\ii(Z): 
h_s(t):=h(s,t)$ and $h^t\in L_\ii(Y): h^t(s):=h(s,t)$, the equality $ess\;sup|h|= 
ess\;sup[ess\;sup|h_s|]= ess\;sup[ess\;sup|h^t|]$. $\tr$ 
%\[
% ess\;sup(h)= ess\;sup[ess\;sup(h_s)]= ess\;sup[ess\;sup(h^t)],
%\]
%где $ess\;sup(h_s)$ рассмотрена как функция от $s\in Y$, а $ess\;sup(h^t)$ -- как функция от $t\in 
%Z$. $\tr$ 

\medskip
We recall that the norm in the injective tensor product $E\ot_{in}F$ of two normed spaces can be 
expressed in terms of an isometric operator from $E\ot_{in} F$ into ${\cal B}(E',F)$, where $E'$ is 
an arbitrary subspace in $E^*$ such that for every $x\in E$ we have $\|x\|=\sup\{|f(x); f\in E', 
\|f\|=1\}$. (For example, $E'$ can be a dense subspace in $E^*$ or a predual to $E$, if such a 
predual exists). See, e.g.,~\cite[pp. 62-63]{clm}) or~\cite[\S 4]{def},~\cite[pp. 45-46]{rya} (and 
also, of course,~\cite{gro}). In the particular case of spaces $L_p(\cd)$, the relevant assertion 
acquires the following form. Denote by $\la f,g\ra$ the classical duality 
$(f,g)\mt\int_Yf(t)g(t)dt$ between the spaces $L_p(Y)$ and $L_q(Y)$, and denote by $\la\la 
u,v\ra\ra$ the duality between the spaces $L_p(Y)\ot^p L_q(Z)$ and $L_q(Y)\ot^q L_q(Z)$, well 
defined on elementary tensor with the help of the equality $\la\la y_1\ot z_1, y_2\ot z_2\ra\ra=\la 
y_1,y_2\ra\la z_1,z_2\ra$. In the following proposition $Y, Z, Y_1,...$ are arbitrary measure 
spaces. 

\medskip
{\bf Proposition 6.2.} {\it 
 (i) There exists an isometric operator ${\cal I}:L_p(Y)\ot_{in}L_q(Z)\to\bb(L_p(Z),L_p(Y))$, (in 
particular, ${\cal I}:{\bf L}L_q(Z)\to\bb(L_p(Z),{\bf L})$), well defined by taking $y\ot z$ to the 
operator, acting by the rule $z'\mt\la z',z\ra y$. 

(ii) If, in addition, $1<p<\ii$, (or, equivalently, $1<q<\ii$), then there exists an isometric 
operator ${\cal J}:[L_p(Y_1)\ot^p L_p(Z_1)] \ot_{in}[L_q(Y_2)\ot^q L_q(Z_2)]\to\bb(L_p(Y_2)\ot^p 
L_p(Z_2),L_p(Y_1)\ot^p L_p(Z_1)$, well defined by taking $u\ot v$ to the operator, acting by the 
rule $v'\mt\la\la v',v\ra\ra u$}. 
 
\smallskip
$\tl$ (i). It is valid because $L_p(Z)$ is the dual (or the predual provided $p=1$) to $L_q(Z)$. 

(ii). It is valid because in the case $1<p<\ii$ the dual space of $(L_q(Y_2)\ot^q L_q(Z_2))$ 
coincides with $L_p(Y_2\times Z_2)$, and the latter space contains $L_p(Y_2)\ot^p L_p(Z_2)$ as a 
dense subspace. $\tr$  

\medskip
{\bf Proposition 6.3.} {\it Let $Y$ and $Z$ be measure spaces, $1<p<\ii$. Then the bilinear 
operator 
 $\rr:L_q(Y)\times L_q(Z)\to L_q(Y)\ot^q L_q(Z)\subset L_q(Y\times Z)$, taking a pair
$(y,z)$ to $y\ot z=y(s)z(t); s\in Y, t\in Z$, is ${\bf L}$--contractive.} 

\smallskip
$\tl$ Consider the bilinear operator ${\cal S}:{\bf L}L_q(Y)\times{\bf L}L_q(Z)=(L_p(X)\ot_{in} 
L_q(Y))\times(L_p(X)\ot_{in} L_q(Z))\to[L_p(X)\ot^p L_p(X)]\ot_{in}[L_q(Y)\ot^q L_q(Z)]$, well 
defined on elementary tensors by the rule $(\xi\ot y,\eta\ot z)\mt(\xi\ot\eta)\ot(y\ot z)$. Take 
$u\in L_p(X)\ot L_q(Y), \\ v\in L_p(X)\ot L_q(Z)$. It is clear that $\rr_\ii=(i_0\ot\id){\cal S}$, 
where  $i_0:L_p(X)\ot^p L_p(X)\to L_p(X)$ is the respective restriction of the isometric 
isomorphism $i$, introduced in Example 3.3, and $\id$ is the identity operator on $L_q(Y)\ot^q 
L_q(Z)$. Since we deal with the injective tensor product, $i_0\ot\id$ is an isometry together with 
$i_0$ and $\id$ (``the injtctive property'', see, e.g.,~\cite[\S 4]{def}). Therefore it is 
sufficient to show that ${\cal S}$ is contractive.  

Consider the diagram 
\[
\xymatrix@C+20pt{(L_p(X)\ot_{in}L_q(Y))\times(L_p(X)\ot_{in}L_q(Z)) \ar[r]^{{\cal S}}\ar[d]_{{\cal I}_1\times {\cal I}_2}
& [L_p(X)\ot^p L_p(X)]\ot_{in}[L_q(Y)\ot^q L_q(Z)] \ar[d]^{{\cal J}} \\
%\bb(L_p(Y),L_p(X))\times\bb(L_p(Z),L_p(X)) \ar[r]^{{\cal T}} & \bb(L_p(Y\times Z), L_p(X)\ot^p L_p(X)) },%\eqno(4)
\bb(L_p(Y),L_p(X))\times\bb(L_p(Z),L_p(X)) \ar[r]^{{\cal T}} & \bb(L_p(Y)\ot^p L^p(Z), L_p(X)\ot^p L_p(X)) },%\eqno(4)
\]
where ${\cal T}$ is the bilinear operator, taking a pair $(A,B)$ to $A\ot B$, and ${\cal I}_k; 
k=1,2$ and ${\cal J}$ are the respective special cases of isometric operators ${\cal I}$ and ${\cal 
J}$ from Proposition 6.2. As it is easy to verify, this diagram is commutative. But by virtue of 
Proposition 6.1 ${\cal T}$ is contractive. Hence, ${\cal S}$ has the same property. $\tr$ 

\medskip
{\bf Theorem 6.4}. {\it Let $Y$ and $Z$ be measure spaces, and $1<p<\ii$. Then we have 
$L_q(Y)\ot_{p{\bf L}}L_q(Z)=L_q(Y)\ot^q L_q(Z)$.} (Recall that the right side of the equality 
denotes the dense subspace in $L_q(Y\times Z)$, consisting of degenerate funnctions). {\it The 
equality means that there exists the ${\bf L}$--isometric isomorphism, well defined by the rule 
$x\ot y\mt x(s)y(t)$ (in other words, it takes $x\ot y$ to the same elementary tensor, only 
considered in $L_q(Y)\ot^q L_q(Z)$).} 

\smallskip
$\tl$ Since all participating spaces are $p$--convex (Proposition 5.4), the Existence Theorem 5.13, 
together with Proposition 6.3, provides the ${\bf L}$--contractive operator $R:L_q(Y)\ot_{p{\bf 
L}}L_q(Z)\to L_q(Y)\ot^q L_q(Z)$, acting as it is indicated in the formulation. We must only show 
that its amplification does not decrease norms. 

Of course, every $L_q(\cd);q<\ii$ contains the dense subspace $L_q^0(\cd)$, consisting of linear 
combinations of characteristic functions of subsets with finite measure, and we have right to 
assume that these subset are pairwise disjoint. 
% and of characteristic functions of single atoms. 
Therefore, because of the estimate (5.7), it is sufficient to show that $R_\ii$ does not decrease 
norms of (finite) sums of elementary tensors of the form $\xi(y\ot z)$, где $\xi\in L_p(X), y\in 
L_q^0(Y), z\in L_q^0(Z)$. 

It is easy to see that every mentioned sum can be represented as $U:=\sum_{k,l}\xi_{k,l}y_k\ot 
z_l$, where $y_k$ and $z_l$ are functions of norm 1, proportional to the characteristic functions 
of pairwise non-intersecting subsets $Y_k\subseteq Y$ and $Z_l\subseteq Z$, respectively.

Look at $R_\ii(U)$. It is, of course, a certain sum $\sum_{k,l}\xi_{k,l}e_{k,l}$, where $e_{k,l}\in 
L_q(Y\times Z)$ are functions of norm 1, proportional to characteristic functions of subsets  
$Y_k\times Z_l\subseteq Y\times Z$; in other words, $e_{k,l}=y_k\ot z_l\in L_q(Y)\ot^q 
L_q(Z)\subset L_q(Y\times Z)$. Therefore, by virtue of Proposition 6.2(i), $\|R_\ii(U)\|$ is the 
norm of the operator $S:L_p(Y\times Z)\to L_p(X)$, taking $h(s,t)$ to $\sum_{k,l}\la h, 
e_{k,l}\ra\xi_{k,l}$. 

Now return to our initial $U$. Obviously, there exist pairwise non-intersecting subsets $X_k^1; 
k=1,...,n$, respectively, $X_l^2; l=1,...,m$ of finite measure.
%, each of them is either non-atomic or is a single atom. 

Denote by $\eta_k$ the functions of norm 1 in ${\bf L}:=L_p(X)$, proportional to the characteristic 
functions of subsets $X_k^1$, and set $u:=\sum_{k=1}^n\eta_ky_k\in{\bf L}L_q(Y)$. By virtue of the 
same Proposition 6.2, $\|u\|$ is the norm of the operator $S_u:L_p(Y)\to L_p(X)$, taking $g$ to 
$\sum_{k=1}^n\la g, y_{k}\ra\eta_k$. Since we obviously have 
$S_u=\sum_{k=1}^nP_{X_k^1}S_{u,k}P_{Y_k}$, where $S_{u,k}:g\mt\la g, y_{k}\ra\eta_k$, we obtain 
that $\|S_u\|=\max\{\|S_{u,k}\|; k=1,...,n\}$ (Proposition 5.2). But we have 
$\|S_{u,k}(g)\|\le\|g\|\|y_{k}\|\|\eta_{k}\|$; therefore $\|S_{u,k}\|\le1$, hence $\|S_u\|\le1$. 
Thus, we have $\|u\|\le1$. Using the same argument, 
%we choose some pairwise disjoint subsets $X_l^2\subset X; l=1,...,m$, where each of them is of the same type as $Z_l$. Then 
we set $v:=\sum_{l=1}^m\zeta_lz_l\in{\bf L}L_q(Z)$ for similarly chosen $\zeta_l$ and we see that 
$\|v\|\le1$. (As a matter of fact, we have $\|u\|=\|v\|=1$, but we do not need it now). 
 
It is clear that $u\di v=\sum_{k,l}(\eta_k\di\zeta_l)y_k\ot z_l$. Therefore it follows from the 
definition of the ${\bf L}$--norm on $L_q(Y)\ot_{p{\bf L}}L_q(Z)$ (see (5.3)) that the theorem will 
be proved, if we shall find an operator $T\in\bb:=\bb(L_p(X))$ such that  
$T(\eta_k\di\zeta_l)=\xi_{k,l}$ (hence, $U=T\cd(u\di v)$) and such that  $\|T\|\le\|S\|$. 

Choose a contractive
%Now observe that there exists an isometric 
operator $i_{k,l}:L_p(X_k^1\times X_l^2)\to L_p(Y_k\times Z_l)$, taking the constant function 
$\eta_k(s)\zeta_l(t); s\in X_k^1, t\in X_l^2$ to a certain constant function $e^*_{k,t}$ of the 
same norm 1. (For example, one can choose $i_{k,l}:=gf$, where $f$ is a functional of norm 1, 
taking the first of mentioned constants to $\equiv 1$, and $g$ is the operator, taking $\equiv 1$ 
to $e^*_{k,t}$). After this, identifying each $L_p(Y_k\times Z_l)$ and $L_p(X_k^1\times X_l^2)$ 
with the respective subspace in $L_p(Y\times Z)$ or, according to the sense, in $L_p(X\times X)$, 
we introduce the operator $D:=S(\sum_{k,l}Q_{k,l}i_{k,l}P_{k,l}):L_p(X\times X)\to L_p(X)$, where 
$P_{k,l}$ is the natural 
   projection of $L_p(X\times X)$ on $L_p(X_k^1\times X_l^2)$, and  $Q_{k,l}$  is the natural
   projection of $L_p(Y\times Z)$ on $L_p(Y_k\times Z_l)$. We see that the operators
    $\sum_{k,l}Q_{k,l}i_{k,l}P_{k,l}$ and $i_{k,l}$ take $\eta_k(s)\zeta_l(t)$ to the same 
constant $e^*_{k,t}$, and the obvious equality $\la e^*_{k,l}, e_{k,l}\ra=1$ implies that 
$S(e^*_{k,l})=\xi_{k,l}$. Therefore $D$ takes $\eta_k(s)\zeta_l(t)$ to $\xi_{k,l}$. 

Further, by virtue of Proposition 5.2 we obtain that $\|\sum Q_{k,l}i_{k,l}P_{k,l}\|\le 
\max_{k,l}\{\|i_{k,l}P_{k,l}\|\}=1$. It follows that $\|D\|\le\|S\|$. 
%the operator $\sum_{k,l}i_{k,l}P_{k,l}$ acts on $\eta_k(s)\zeta_l(t)$ exactly as $i_{k,l}$, and it follows from 
%the obvious equality $\la e^*_{k,l}, e_{k,l}\ra=1$ that $S(e^*_{k,l})=\xi_{k,l}$. Therefore $D$ 
%takes $\eta_k(s)\zeta_l(t)$ to $\xi_{k,l}$. 

%It is clear that $i_{k,l}P_{k,l}=Q_{k,l}i_{k,l}P_{k,l}$, where $Q_{k,l}$ is the natural projection 
%of $L_p(Y\times Z)$ onto $L_p(Y_k\times Z_l)$. From this, by virtue of Proposition 5.2, we obtain 
%that $\|\sum_{k,l}i_{k,l}P_{k,l}\|\le \max\{\|i_{k,l}P_{k,l}\|\}=1$. Consequently, we have 
%$\|D\|\le\|S\|$. 
Finally, let us recall the isometric isomorphism $ i:L_p(X\times X)\to L_p(X)$ (see Example 3.3). 
This map, in particular, takes every function of the form $\eta(s)\zeta(t)$ (identified with 
$\eta\ot\zeta\in L_p(X)\ot^p L_p(X)$) to $\eta\di\zeta$. From this we see that the operator 
$T:=Di^{-1}:L_p(X)\to L_p(X)$ is exactly what we need. $\tr$ 
 
\medskip 
{\bf Remark 6.5.} Throughout all paper, we did not assume that our normed spaces are complete. 
However, principal notions and facts have, as a rule, ``complete'' versions. ${\bf L}$--space is 
called complete (or Banach), if its underlying normed space is complete. As in the ``classical'' 
context, every ${\bf L}$--space has a completion; its definition and properties, as well as the 
existence theorem repeat with obvious modifications what was said in~\cite[Ch. 4]{heb2} for the 
case of operator spaces.  (See also similar argument in~\cite{heh} for the case, where 
 ${\bf L}$ is a Hilbert space). Moreover, both tensor products, introduced above, have their ``Banach'' 
versions; one must only to consider in Definition 4.2, in the capacity of $\mho$, the class of 
complete ${\bf L}$-spaces, whereas in Definition 5.5 one must consider the class of complete 
$p$--convex ${\bf L}$--spaces. Then, in particular, we obtain the following version of Theorem 6.5: 
if ${\bf L}:=L_p(X)$ and $1<p<\ii$, then {\it completed} $p$--convex tensor product of ${\bf 
L}$-spaces $L_q(Y)$ and $L_q(Z)$ is  $L_q(Y\times Z)$. 

%\newpage

\ed